\documentclass{amsart}
\usepackage[final]{epsfig}
\usepackage{graphics}
\usepackage{float}
\usepackage{amsmath}
\usepackage{amsfonts}
\usepackage{latexsym}
\usepackage{amssymb}
\usepackage{graphicx}
\usepackage{url}
\usepackage{epstopdf}
\usepackage{verbatim}
\usepackage[margin=1in]{geometry}
\usepackage{amsthm}
\usepackage{enumitem}
\usepackage{tikz}
\usetikzlibrary{shapes,arrows}

\newtheorem{lemma}{Lemma}

\newtheorem{theorem}{Theorem}

\newtheorem*{definition*}{Definition}

\newtheorem*{theorem*}{Theorem}

\newcommand{\proofend}{$\Box$\bigskip}

\newcommand{\R}{{\mathbb R}}

\newcommand{\s}{{\mathbb S}}
\newcommand{\B}{{\mathbb B}}

\def\proof{\paragraph{Proof.}}

\begin{document}

\title{A Simple Proof of Cauchy's Surface Area Formula}

\author{Emmanuel Tsukerman}
\address{Department of Mathematics, University of California,
Berkeley, CA 94720-3840}
\email{e.tsukerman@berkeley.edu}
\author{Ellen Veomett}
\address{Department of Mathematics and Computer Science, Saint Mary's College of California, Moraga, CA, 94575}
\email{erv2@stmarys-ca.edu}

\date{\today}

\begin{abstract}
We give a short and simple proof of Cauchy's surface area formula, which states that the average area of a projection of a convex body is equal to its surface area up to a multiplicative constant in the dimension. 
\end{abstract}

\maketitle
\setcounter{tocdepth}{1}

\section{Introduction}

Cauchy's surface area formula is a beautiful and well-known result in integral geometry.  It states that the average area of the projections of a convex body is equal to the surface area of the body, up to a multiplicative constant in the dimension.  That is, let $\mu_{n-1}$ be $(n-1)$-dimensional Lebesgue measure on $\R^{n-1}$, and  $\s^{n-1}$ be the unit sphere in $\R^n$.  For a convex body $K \subset \R^n$ and 
$u \in \s^{n-1}$, we denote by $K|u^{\perp}$ the projection of $K$ onto the $(n-1)$-dimensional subspace of $\R^n$ perpendicular to $u$.  Then we have:

\begin{theorem}[Cauchy's surface area formula]\label{CauchyTheorem}
Suppose $K \subset \R^n$ is a convex body.  Then
\begin{equation}
\text{S}(K) =  \frac{1}{\mu_{n-1}\left(\B^{n-1}\right)} \int_{\s^{n-1}} \mu_{n-1}\left(K|u^\perp \right) du.
\end{equation}
where $S(K)$ denotes the volume of the surface of $K$.
\end{theorem}

We note that $S(K)$ is succinctly defined as follows:
\begin{equation}\label{SurfaceArea}
S(K) = \lim_{\epsilon \to 0^+} \frac{\mu_n\left(K+\epsilon \B_n\right)-\mu_n(K)}{\epsilon}
\end{equation}
Where $\B_n$ is the unit ball in $\R^n$ and $K+\epsilon B_n$ indicates the Minkowski sum of $K$ and $\epsilon B_n$:
\begin{equation*}
K+\epsilon \B_n = \{x+\epsilon v: x \in K, v \in \B_n\}.
\end{equation*}

Cauchy's surface area formula was first proved by Cauchy for $n=2, 3$  in 1841 and 1850 \cite{Cauchy1}, \cite{Cauchy2}.  The general formula was proved in work by Kubota \cite{Kubota1}, and also Minkowski \cite{Minkowski1}  and Bonnesen \cite{Bonnesen1}.  Klain and Rota give a lovely proof in Chapter 5 of \cite{MR1608265}, which uses calculus-like arguments and approximations of convex bodies by polytopes.  (We draw much of our notation from the same text).  In this short note we give an alternative proof, which we believe is new.  This new proof has  the benefit of the fact that the projection is calculated using a Minkowski sum of sets, just as the surface area formula is calculated using a Minkowski sum of sets as in equation \eqref{SurfaceArea}.  Thus, the two sides of Theorem \ref{CauchyTheorem} are seen to be calculated using Minkowski sums.  This technique allows  a straightforward extension to intrinsic moment vectors of a convex body (which we will define in Section \ref{MomentVectors}).

\section{New Proof for Cauchy's Surface Area Formula}\label{NewProof}

First we give our definition, which uses the Minkowski sum of sets as discussed in the Introduction.

\begin{definition*}
Let $u, P \subset \R^n$. Define
\begin{equation}\label{deriv}
D_u  (\mu_n)(P)=\lim_{\epsilon \rightarrow 0^{+}} \frac{\mu_n(P+\epsilon u)-\mu_n(P)}{\epsilon}.
\end{equation}
\end{definition*}

The special case $D_{\mathbb{B}^n}  (\mu_n)(P)$ gives the surface volume:
\begin{equation}\label{minkcontent}
D_{\mathbb{B}^n}  (\mu_n)(P) =S(P)
\end{equation}

The key to our proof is Minkowski's theorem on mixed volumes, which can be found, for example, in  Chapter 5 of Schneider's text \cite{MR3155183}.  This theorem says that the volume of a Minkowski sum of convex bodies can be written as a polynomial in the coefficients of the Minkowski sum, where the coefficients of the polynomial depend only on the convex bodies.  Specifically:  

\begin{theorem}\label{MixedVolumeThm}
Suppose $K_1, K_2, \dots, K_m$ are convex bodies in $\R^n$.  Then
\begin{equation*}
\mu_n\left(\lambda_1K_1+\lambda_2 K_2+ \cdots + \lambda_m K_m\right) = \sum \lambda_{i_1}\lambda_{i_2} \cdots \lambda_{i_n} V\left(K_{i_1}, K_{i_2}, \dots, K_{i_n}\right)
\end{equation*}
where the sum on the left-hand side is the Minkowski sum, and the sum on the right-hand side is over all multisets of size $n$ whose elements are in the set $\{1, 2, \dots, m\}$.  The  functions $V$ are nonnegative, symmetric, and depend only on the convex bodies $K_{i_1}, K_{i_2}, \dots, K_{i_n}$.
\end{theorem}

From Minkowski's Theorem it follows that if $P$ and $u$ are convex, then $D_u  (\mu_n)(P)$ is linear in $u$:

\begin{lemma}\label{LinearityLemma}
Suppose $P, u, v \subset \R^n$ with $P, u$, and $v$ convex bodies, and suppose $\alpha, \beta \in \R$.  Then
\begin{equation}
D_{\alpha u+\beta v} (\mu_n)(P) = \alpha D_u (\mu_n)(P)+\beta D_v (\mu_n)(P)
\end{equation}
\end{lemma}

\proof
By definition we have
\begin{equation*}
D_{\alpha u+\beta v}(\mu_n)(P) = \lim_{\epsilon \to 0^+} \frac{\mu_n\left(P+\epsilon \left(\alpha u+\beta v\right)\right) - \mu_n(P)}{\epsilon} =  \lim_{\epsilon \to 0^+} \frac{\mu_n\left(P+\epsilon \alpha u+ \epsilon \beta v\right) - \mu_n(P)}{\epsilon} 
\end{equation*}
From Theorem \ref{MixedVolumeThm}, we can see that
\begin{equation*}
D_{\alpha u+\beta v}(\mu_n)(P) = \alpha V(\underbrace{P, P, \dots, P}_{n-1 \text{ times}}, u) + \beta V(\underbrace{P, P, \dots, P}_{n-1 \text{ times}}, v) 
\end{equation*}
where $V$ represents the function in the statement of Theorem \ref{MixedVolumeThm}.  Similarly, one can easily see that
\begin{align*}
D_u(\mu_n)(P) &=  V(\underbrace{P, P, \dots, P}_{n-1 \text{ times}}, u) \\
 D_v(\mu_n)(P) &=  V(\underbrace{P, P, \dots, P}_{n-1 \text{ times}}, v), 
 \end{align*}
and our Lemma is proved.
\proofend

 For the following Lemma, we prove that the derivative we've defined in equation \eqref{deriv}  calculates the projection of a convex body in the case where $u$ is a segment of    length 1.  We use the same notation as in the Introduction.

\begin{lemma}\label{DerivativeLemma}
Let $u$ be a segment of length $1$ and $K \subset \R^n$ convex. Then
\begin{equation}\label{projeq}
\mu_{n-1} \left(K|u^\perp\right) = D_u(\mu_n)( K).
\end{equation}
\end{lemma}

\proof
Let $L_u$ be the set of lines parallel to $u$.  Note that each $l \in L_u$ corresponds uniquely to a single point in $l^\perp$, so that $L_u$ is isomorphic to $\R^{n-1}$ (and thus we can define the measure $\mu_{n-1}$ on $L_u$).  Then
\[
\mu_n(K)=\int_{l \in L_u} \mu_1(l \cap K) \, d \mu_{n-1}.
\]
Also,
\[
\mu_{n-1}\left(K|u^\perp\right) =\int_{\substack{ l \cap K \neq \emptyset \\ l \in L_u }} 1 \, d \mu_{n-1}.
\]
For $\epsilon >0$, we have
\[
\mu_n\left(K+\epsilon u\right)-\mu_n\left(K\right)=\int_{l \in L_u}\left( \mu_1\left(l \cap (K+\epsilon u)\right)-\mu_1\left(l \cap K\right) \right)d\mu_{n-1}.
\]
Convexity implies that
\[
\mu_1\left(l \cap (K+\epsilon u)\right)-\mu_1\left(l \cap K\right) = \epsilon
\]
whenever $l$ intersects $K$, and zero otherwise. Therefore the last integral is equal to
\[
\int_{\substack{ l \cap K \neq \emptyset \\ l \in L_u }} \epsilon \, d \mu_{n-1}=\epsilon\mu_{n-1} \left(K|u^\perp\right) 
\]
and hence,
\begin{equation*}
D_u (\mu_n)(K) = \lim_{\epsilon \to 0^+}  \frac{\mu_n(K+\epsilon u)-\mu_n(K)}{\epsilon} = \lim_{\epsilon \to 0^+} \frac{\epsilon\mu_{n-1} \left(K|u^\perp\right) }{\epsilon} = \mu_{n-1} \left(K|u^\perp\right). 
\end{equation*}
\proofend

We can now give a short proof of Cauchy's surface area formula:

\proof (of Theorem \ref{CauchyTheorem})
 Using Lemmas \ref{LinearityLemma},\ref{DerivativeLemma}, and equation \eqref{minkcontent}, we have: 
\begin{equation*}
 \int_{\s^{n-1}} \mu_{n-1}\left(K|u^\perp \right) du
= \int_{\s^{n-1}} D_{[0,u]} (\mu_n)( K) du
=  D_{\int_{\s^{n-1}}[0,u]du} (\mu_n)( K)
=    D_{c(n) \B^n}  (\mu_n)(K) 
=  c(n)  S(K). 
\end{equation*}

Specializing $K$ to be the unit ball $\B^n$ yields the result.
\proofend

\section{Application of This Technique to Moment Vectors}\label{MomentVectors}

The proof can be extended readily to intrinstic moment vectors of a convex body \cite[Sec. 5.4]{MR3155183}. We introduce the following notation, modeled on the notation from  \cite[Sec. 5.4]{MR3155183}.  Let $H^r$ denote the $r$-dimensional Hausdorff measure on $\R^n$.   For an $r$ dimensional set $K \subset \R^n$, let 
\begin{equation*}
z_{r+1}(K)=\int_K x \, dH^r(x)
\end{equation*}
 denote the moment vector of the set $K$. It is related to the familiar centroid $p_0(K)$ via $z_{r+1}(K)=H^r(K)p_0(K)$. We will abuse notation and denote the boundary of $K$ by $S(K)$, when no confusion with the surface area can occur. For a segment $u$, let 
\[
K_u= \lim_{\epsilon \rightarrow 0} S(K) \cap (K+\epsilon u).
\]
Intuitively, this is the portion of the boundary that is seen by a projection in direction $u$.

\begin{figure}[H]
\centering
\includegraphics[width=2in]{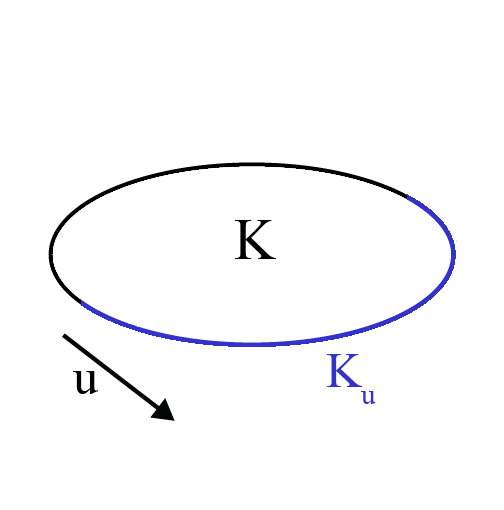}
\caption{The subset $K_u$ of the boundary $\partial K$ of $K$.}
\end{figure}

Let $K$ be an $n$-dimensional subset of $\R^n$.  The following Theorem now has a very similar flavor to that of the Cauchy surface area formula:
\begin{theorem}\label{MomentTheorem}
Suppose $K \subset \R^n$ is a convex body.  Then
\[
z_n(S (K))=\frac{1}{\mu_{n-1}(\B_{n-1})} \int_{\s^{n-1}} z_{n}(K_u) du .
\]
\end{theorem}

Before we prove this, we need a definition and a few observations.  If $K$ is an  $n$-dimensional subset of $\R^n$ and $u$ any subset of $\R^n$, we now define
\begin{equation}\label{MomentDerivative}
D_u (z_{n+1})(K) = \lim_{\epsilon \to 0^+} \frac{z_{n+1}(K+\epsilon u) - z_{n+1}(K)}{\epsilon}.
\end{equation}
We note that, as one can find in  \cite[Sec. 5.4]{MR3155183}, there is a Theorem parallel to Theorem \ref{MixedVolumeThm} but for the moment vector.  Thus, just as we showed in Lemma \ref{LinearityLemma}, the operator $D_u$ defined in equation \eqref{MomentDerivative} is also linear in $u$ (if $K, u, $ and $v$ are convex just as in Lemma \ref{LinearityLemma}).  We can now proceed with our proof  of Theorem \ref{MomentTheorem}:

\proof (Proof of Theorem \ref{MomentTheorem})
It is not hard to see from the definition of the moment vector that 
\[
z_{n+1}(K+\epsilon u)=z_{n+1}(K)+z_{n+1}((K+\epsilon u) \setminus K)
\]

Therefore if $u \in \s^{n-1}$ and as above we let $[0,u]$ denote the segment from the origin to $u$, we have
\[
D_{[0,u]}(z_{n+1})(K)=\lim_{\epsilon \rightarrow 0^+} \frac{z_{n+1}((K+\epsilon u) \setminus K)}{\epsilon}=z_n(K_u).
\]
which gives 
\[
\int_{\s^{n-1}} D_{[0,u]}(z_{n+1})(K) du=\int_{\s^{n-1}} z_n(K_u) du.
\]
Similarly, we note that
\begin{equation*}
z_{n+1}(K+\epsilon B^n) = z_{n+1}(K)+z_{n+1}((K+\epsilon \B_n) \backslash K)
\end{equation*}
so that
\begin{equation}\label{DerivativeMomentBall}
D_{\B^n} (z_{n+1})(K) = \lim_{\epsilon \rightarrow 0^+} \frac{z_{n+1}((K+ \epsilon \B_n) \backslash K}{\epsilon} = z_n(S(K)).
\end{equation}

Note that $K$ and $[0,u]$ are convex bodies so that we can use the linearity of the $D$ operator defined in \eqref{MomentDerivative}.  Using this linearity, our calculations from the proof of Theorem \ref{CauchyTheorem}, and equation \eqref{DerivativeMomentBall}, we now have
\begin{multline*}
\int_{\s^{n-1}} D_{[0,u]}(z_{n+1})(K)  \,du=D_{\int_{ \s^{n-1}}  [0,u] du} (z_{n+1})(K)
= \mu_{n-1}(\B^{n-1})D_{\B^n}(z_{n+1})(K) = \mu_{n-1}(\B^{n-1})z_n(S(K)).
\end{multline*}
Thus the Theorem is proved.
\proofend

\bigskip
{\bf Acknowledgments}.  
The authors would like to thank Alexander Barvinok for his comments and Sergei Tabachnikov for his advice and encouragement.

This material is based upon work supported by the National Science Foundation Graduate Research Fellowship under Grant No. DGE 1106400. Any opinion, findings, and conclusions or recommendations expressed in this material are those of the authors(s) and do not necessarily reflect the views of the National Science Foundation.

\bibliographystyle{alpha}
\bibliography{AvgProj}

\end{document}